\documentclass{amsart}
\usepackage{amsmath}
\usepackage{amsfonts}
\usepackage{amssymb}
\usepackage[latin2]{inputenc}
\usepackage{graphicx}

\DeclareMathOperator{\co}{\mathbb C}
\DeclareMathOperator{\K}{K}
\DeclareMathOperator{\KX}{K_{1}}
\DeclareMathOperator{\KY}{K_{\bar{1}}}
\DeclareMathOperator{\KXX}{K_{11}}
\DeclareMathOperator{\KYY}{K_{\bar{1}\bar{1}}}
\DeclareMathOperator{\KXY}{K_{1\bar{1}}}
\DeclareMathOperator{\KXXY}{K_{11\bar{1}}}
\DeclareMathOperator{\KXYY}{K_{1\bar{1}\bar{1}}}
\DeclareMathOperator{\KXXYY}{K_{11\bar{1}\bar{1}}}
\DeclareMathOperator{\mianown}{\left(-\frac{\KY \KX}{\K^2}+\frac{\KXY}{\K}\right)}
\DeclareMathOperator{\mian}{S}
\newtheorem{theorem}{Theorem}[section]
\newtheorem{lemma}[theorem]{Lemma}
\newtheorem{observation}{Observation}
\theoremstyle{definition}

\numberwithin{equation}{section}
\begin{document}
\title[the holomorphic sectional curvature]{An example for the holomorphic sectional curvature of the Bergman metric}
\author{Żywomir Dinew}
\address{Jagiellonian University, Institute of Mathematics,
Łojasiewicza 6,
30-348 Kraków,
Poland  }
\email{Zywomir.Dinew@im.uj.edu.pl}
\thanks{The author was supported  by the Polish Ministry of Science and Higher Education Grant N N201 271235.}
\begin{abstract}
In this paper we study the behaviour of the holomorphic sectional curvature (or Gaussian curvature) of the Bergman metric of planar annuli. The results are then utilized to construct a domain for which the curvature is divergent at one of its boundary points and moreover the limes superior is the maximal possible for the Bergman metric (2), whereas the limes inferior is $-\infty$.  
\end{abstract}
\maketitle
\begin{section}{Introduction}
Recall, that the holomorphic sectional curvature of the Bergman metric of a bounded pseudoconvex domain $U\subset\mathbb C^{n}$, at the point $z\in U$ in direction $X\in\co^{n}$ is defined as follows:

\begin{equation}\label{krzywizna}
R_{U}(z,X):=\left(\sum_{p,q=1}^{n}g_{p\overline{q}} X_{p}X_{\overline{q}}\right)^{-2}\sum_{i,j,k,l=1}^{n} R_{\overline{i}jk\overline{l}}\overline{X_i}X_jX_k\overline{X_l},\end{equation}
here 
$$R_{\overline{i}jk\overline{l}}:=-\frac{\partial^{2}g_{j\overline{i}}}{\partial z_{k}\partial\overline{z}_{l}}+\sum_{r,s=1}^{n}g^{\overline{r}s}\frac{\partial g_{j\overline{r}}}{\partial z_{k}}\frac{\partial g_{s\overline{i}}}{\partial \overline{z}_{l}},$$
where $g^{\overline{r}s}$ stands for $r,s$-th element of the inverse matrix of $g_{p\overline{q}}$. The term in brackets in the definition of $R_{U}$ is introduced for the sake of normalization. Finally
 $g_{p\overline{q}}$ stands for $\frac{\partial ^2}{\partial z_{p}\partial\overline{z}_{q}}\log \K_{U}(z,z),$
where $\K_{U}(z,z)$ is the Bergman kernel (on the diagonal) of the domain $U$.

One can show that:
\begin{equation}\label{r}R_{U}(z_{0},X):=2-\frac{J_{1,U}(z_{0};X)^2}{J_{0,U}(z_{0};X)J_{2,U}(z_{0};X)},\end{equation}
where
\begin{equation}\label{j0} J_{0,U}(z_{0};X):=sup\{|f(z_{0})|^2: f(z)\in\mathcal O(U)\cap L^2(U), \int_{U} |f|^2\leq1\}
\end{equation}
 $$ J_{1,U}(z_{0};X):=sup\Big\{\Big|\sum_{j=1}^{n}\frac{\partial f(z_{0})}{\partial z_{j}}X_{j}\Big|^2: f(z)\in\mathcal O(U)\cap L^2(U),$$
\begin{equation}\label{j1}  \int_{U} |f|^2\leq1, f(z_{0})=0\Big\}
\end{equation}
 $$J_{2,U}(z_{0};X):=sup\Big\{\Big|\sum_{i,j=1}^{n}\frac{\partial^2 f(z_{0})}{\partial z_{j}\partial z_{i}}X_{j}X_{i}\Big|^2: f(z)\in\mathcal O(U)\cap L^2(U),$$ 
\begin{equation}\label{j2} \int_{U} |f|^2\leq 1, f(z_{0})=0, \frac{\partial f(z_{0})}{\partial z_{j}}=0,j=1..n\Big\}
\end{equation}

We see that $J_{0,U}(z_{0};X)=\K_{U}(z_{0},z_{0})$, which is independent of $X$, and that the holomorphic sectional curvature of the Bergman metric is invariant under biholomorphic mappings. 

From \ref{r} it follows immediately that
\begin{equation}\label{leb}R_{U}(z,X)<2,z\in U.
\end{equation}
 This was known already by Bergman. It was shown by Lebed (see \cite{MR0296346}) that when $n\geq2$ this estimate is optimal in the following weak sense: For each $\varepsilon>0$ there exists a domain $U_{\varepsilon}$, for which there exist $z\in U_{\varepsilon}$ and $X\in \mathbb C^{n}$, such that $R_{U_{\varepsilon}}(z,X)>2-\varepsilon$. In a very recent paper  Chen and Lee (\cite{chenlee}) have shown that the estimate is optimal in the strong sense, i.e., there exists a domain $U$ and $z_{0}\in\partial U$, such that $\lim_{\nu \to \infty}R_{U}(z_{\nu},X(z_{\nu}))=2$ for suitably chosen $z_{\nu}\in U,z_{\nu}\to z_{0},X(z_{\nu})$. The question of the existence of a lower bound is also answered (in the negative) in higher dimensions in the paper of Herbort (see \cite{MR2318508}, where even an example with smooth boundary is provided ). 

In dimension $1$ the formula \eqref{krzywizna} becomes 
\begin{equation}\label{onedim}
R_{U}(z,X)=\frac{-g_{\overline{1}11\overline{1}}+\frac{g_{\overline{1}1\overline{1}}g_{1\overline{1}1}}{g_{1\overline{1}}}}{(g_{1\overline{1}})^2}=\frac{-(\log{g_{1\overline{1}}})_{1\overline{1} }}{g_{1\overline{1}}},\end{equation}
which is independent of $X$ and therefore we will use the shorter form $R_{U}(z)$. In fact this is exactly the Gaussian curvature of $g$.

Little is known about the holomorphic sectional curvature of the Bergman metric in dimension one. This is mainly because one cannot compute the Bergman kernel of most of the  planar domains explicitly. The first nontrivial (i.e., not biholomorphic to the unit disc) domain for which one can say more is the circular annulus.

The Bergman kernel of the annulus $P_{r}=\{z\in\co: r<|z|<1\}$ is
$$\K_{P_{r}}(z,z)=-\frac{1}{\pi  |z|^2 \log(r^2)}+\pi^{-1} \sum _{j=0}^{\infty } \left(\frac{r^{2+2 j}}{\left(-r^{2+2 j}+|z|^2\right)^2}+\frac{r^{2 j}}{\left(1-r^{2
j} |z|^2\right)^2}\right)$$
(see e.g. \cite{MR1242120}). This will be denoted as $\K$ for short.

Because the annulus has smooth boundary, it follows by a theorem of Klembeck (see \cite{MR0463506}), that
$$\lim_{P_{r}\ni z\to\partial P_{r}}R_{P_{r}}(z)=-1.$$ 
When $r\to 0$, the domains $P_{r}$ will exhaust the punctured unit disc and therefore one can expect the corresponding holomorphic sectional curvatures of the Bergman metric to be convergent to the the curvature of the punctured disc, which is the same as for the whole disc (the constant $-1$). (In the case of the Bergman kernel this is the assertion of the Ramadanov's theorem). This is indeed the case, however the convergence is only locally uniform. Moreover, numerical experiments of the author have shown that the global maximum of the holomorphic sectional curvature of the Bergman metric becomes closer to $2$, the smaller $r$ is and the global minimum tends to be unbounded.

The figures below present the behaviour of the curvature, when restricted to the line segment $(r,1)\subset \mathbb R$, for different choices of $r$.(Figures $3$ and $4$ are for the same $r=0.001$, however figure $4$ is scaled in order to stress on the global maximum).

\begin{figure} [tbhp]
\centering
\includegraphics[width=2.31 in]{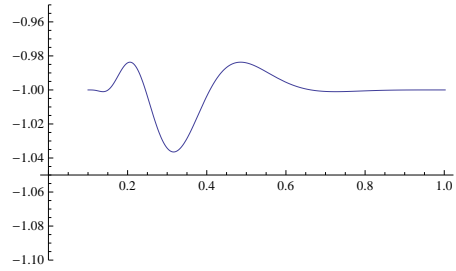}
\caption{The curvature of $P_{0.1}$ restricted to (0.1,1)}
\label{r=0.1}
\end{figure}
\begin{figure}[tbph]
\centering
\includegraphics[width=2.31 in]{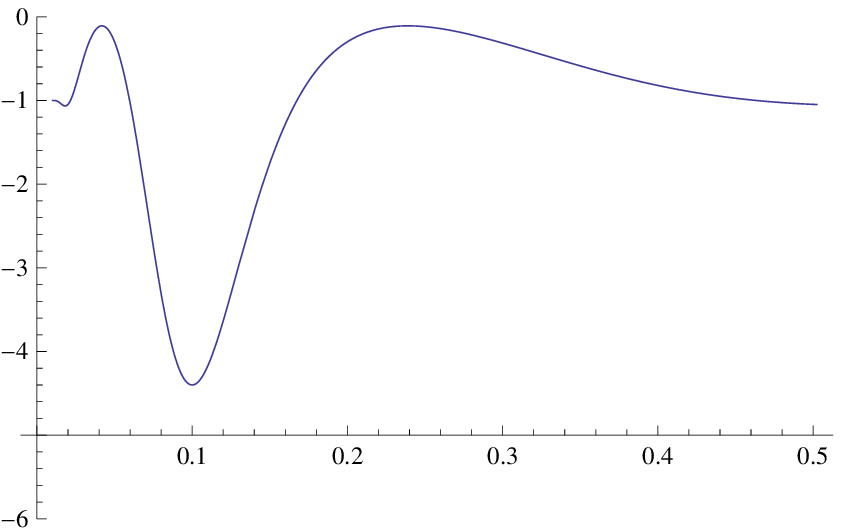}
\caption{The curvature of $P_{0.01}$ restricted to (0.01,1)}
\label{r=0.01}
\end{figure}
\begin{figure}[tbhp]
\centering
\includegraphics[width=2.31 in]{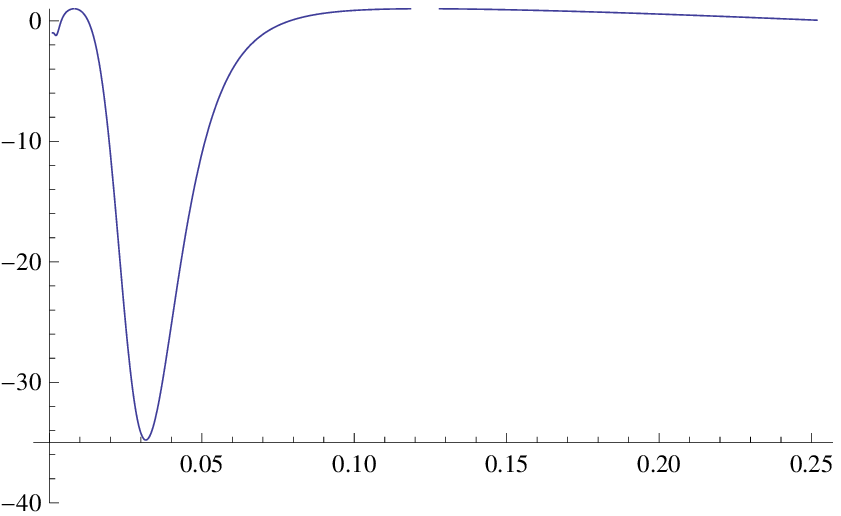}
\caption{The curvature of $P_{0.001}$ restricted to (0.001,1)}
\label{r=0.001}
\end{figure}
\newpage
\begin{figure}[tbph]
\centering
\includegraphics[width=2.31 in]{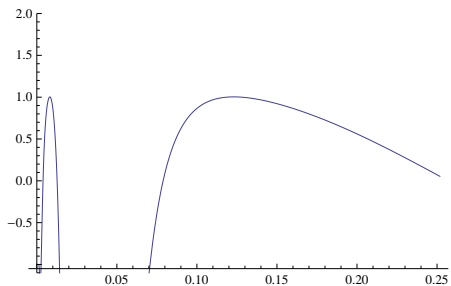}
\caption{The curvature of $P_{0.001}$ restricted to (0.001,1)}
\label{r=0.001scaled}
\end{figure}

 We will confirm these numerical experiments by proving analitically the following result
\begin{theorem}\label{glownetw} For the circular annulus $P_{r}$ one has

$$\lim_{r\to 0} R_{P_{r}}(\sqrt{r})=-\infty$$
$$\lim_{r\to 0} R_{P_{r}}(r^{\frac{3}{10}})=2$$
 
\end{theorem}

Section $1$ is entirely devoted to the very technical but rigorous proof of this theorem. It is of course desirable to find an easier proof. 

This shows that Lebed's result can be extended to $n=1$ and that one cannot find a universal constant to be a lower bound for any planar domain.

In section $2$  the result of section $1$ is utilized to construct a planar domain for which the holomorphic sectional curvature of the Bergman metric is divergent at one of its boundary points. Namely one has

\begin{theorem} There exists a bounded planar domain $\Omega$ and a point $\zeta\in\partial\Omega$ such that
$$\limsup_{\Omega\ni z\to \zeta}R_{\Omega}(z)=2$$  and $$\liminf_{\Omega\ni z\to \zeta}R_{\Omega}(z)=-\infty.$$ 
\end{theorem}
This is done by using known (in the style of \cite{MR978601}, or \cite{MR2002015}) localization technique, which heavily depends on the geometry of the domain.

This proves that the constant $2$ is optimal in the strong sense even in dimension $1$ and that planar domains with unbounded holomorphic sectional curvature of the Bergman metric do exist.

In the whole paper $G_{\Omega}(z,z_{0})$ will stand for the (pluricomplex) Green function i.e.,
$$G_{\Omega}(z,z_{0}):=sup\{u(z): u\in PSH(\Omega), u\leq0 , \limsup_{w\to z_{0}} u(w)-\log|w-z_{0}|<\infty\}$$
Note that in dimension $1$ it is customary to call  $-G(z,z_{0})$ ``Green function'', however we will stick to the above definition (in dimension $1$ one has to change $PSH$ to $SH$ and $\limsup$ is in fact $\lim$).
\end{section}
\begin{section}{Behaviour of the curvature of circular annuli}
Let $S:=\mianown$. Using \eqref{onedim} we expand explicitly $R_{P_{r}}$ by means of consecutive derivatives of the Bergman kernel $\K$ as

$R_{P_{r}}=\sum_{j=1}^{24}A_{j}$ , where
\begin{align*}
A_{1}&=\frac{4 \KY^3 \KX^3}{\K^6 \mian^3},\ &A_{2}&=-\frac{2\KY \KYY\KX^3}{K^5 \mian^3},\
&A_3&=-\frac{8\KY^2 \KX^2 \KXY}{\K^5 \mian^3},\\
 A_4&=\frac{2
\KYY \KX^2 \KXY}{\K^4 \mian^3},\ &A_5&=\frac{4
\KY\KX\KXY^2}{\K^4 \mian^3},\ &A_6&=\frac{6
\KY^2\KX^2}{\K^4 \mian^2},\\
A_7&=-\frac{2
\KYY\KX^2}{\K^3 \mian^2},\ &A_8&=-\frac{8 \KY
\KX\KXY}{\K^3 \mian^2},\ &A_9&=\frac{2 \KXY^2}{\K^2
\mian^2},\\
A_{10}&=\frac{2 \KY\KX^2\KXYY}{\K^4
\mian^3},\  &A_{11}&=-\frac{2 \KX\KXY\KXYY}{\K^3
\mian^3},\ &A_{12}&=\frac{2 \KX \KXYY}{\K^2 \mian^2}\\
A_{13}&=-\frac{2 \KY^3 \KX \KXX}{\K^5 \mian^3},\  &A_{14}&=\frac{\KY\KYY\KX\KXX}{\K^4 \mian^3},\ &A_{15}&=\frac{2 \KY^2\KXY\KXX}{\K^4 \mian^3},\\
 A_{16}&=-\frac{2 \KY^2\KXX}{\K^3 \mian^2},\ &A_{17}&=\frac{\KXX
\KYY}{\K^2 \mian^2},\ &A_{18}&=-\frac{\KY\KXYY\KXX}{\K^3 \mian^3},\\
A_{19}&=\frac{2 \KY^2 \KX
\KXXY}{\K^4 \mian^3},\  &A_{20}&=-\frac{\KYY \KX \KXXY}{\K^3 \mian^3},\ &A_{21}&=-\frac{2 \KY \KXY
\KXXY}{\K^3 \mian^3},\\
 A_{22}&=\frac{2 \KY \KXXY}{\K^2 \mian^2},\ &A_{23}&=\frac{\KXYY \KXXY}{\K^2 \mian^3},\  &A_{24}&=-\frac{\KXXYY}{\K\mian^2}.
\end{align*}

By direct computation one obtains
$$\K=-\frac{1}{\pi  |z|^2 \log(r^2)}+\pi^{-1} \sum _{j=0}^{\infty } \left(\frac{r^{2+2 j}}{\left(-r^{2+2 j}+|z|^2\right)^2}+\frac{r^{2 j}}{\left(1-r^{2
j} |z|^2\right)^2}\right)$$
$$\KX=\frac{1}{\pi  z|z|^2 \log(r^2)}+\pi^{-1} \sum _{j=0}^{\infty } \left(-\frac{2r^{2+2 j}\overline{z}}{\left(-r^{2+2 j}+|z|^2\right)^3}+\frac{2r^{4j}\overline{z}}{\left(1-r^{2
j} |z|^2\right)^3}\right)$$
$$\KY=\overline{\KX}$$
$$\KXX=-\frac{2}{\pi  z^2|z|^2 \log(r^2)}+\pi ^{-1}\sum _{j=0}^{\infty }\left( 
\frac{6r^{2+2j}\overline{z}^2}{\left(-r^{2+2 j}+|z|^2\right)^4}+
\frac{6r^{6j}\overline{z}^2}{\left(1-r^{2j}|z|^2\right)^4}\right)$$
$$\KYY=\overline{\KXX}$$
$$\KXY=-\frac{1}{\pi  |z|^4 \log(r^2)}+\pi ^{-1}\sum _{j=0}^{\infty } \Big(
\frac{6r^{2+2j}|z|^2}{\left(-r^{2+2 j}+|z|^2\right)^4}-
\frac{2r^{2+2j}}{\left(-r^{2+2j}+|z|^2\right)^3}+$$
$$\frac{6r^{6j}|z|^2}{\left(1-r^{2j}|z|^2\right)^4}+
\frac{2r^{4j}}{\left(1-r^{2j}|z|^2\right)^3}\Big)$$
$$\KXXY=\frac{2}{\pi  z|z|^4 \log(r^2)}+\pi ^{-1}\sum _{j=0}^{\infty } \Big(
-\frac{24r^{2+2j}|z|^2\overline{z}}{\left(-r^{2+2 j}+|z|^2\right)^5}+
\frac{12r^{2+2j}\overline{z}}{\left(-r^{2+2j}+|z|^2\right)^4}+$$
$$
\frac{24r^{8j}|z|^2\overline{z}}{\left(1-r^{2j}|z|^2\right)^5}+
\frac{12r^{6j}\overline{z}}{\left(1-r^{2j}|z|^2\right)^4}\Big)$$
$$\KXYY=\overline{\KXXY}$$
$$\KXXYY=-\frac{4}{\pi  |z|^6 \log(r^2)}+\pi ^{-1}\sum _{j=0}^{\infty }\Big[6r^{2+2j} \Big(
\frac{20|z|^4}{\left(-r^{2+2 j}+|z|^2\right)^6}-$$$$
\frac{16|z|^2}{\left(-r^{2+2 j}+|z|^2\right)^5}+ 
\frac{2}{\left(-r^{2+2 j}+|z|^2\right)^4}\Big) +6r^{6j}\Big(
\frac{20r^{4j}|z|^4}{\left(1-r^{2j}|z|^2\right)^6}+
\frac{16r^{2j}|z|^2}{\left(1-r^{2j}|z|^2\right)^5}+$$$$
\frac{2}{\left(1-r^{2j}|z|^2\right)^4}\Big)\Big]$$

 For the special choice $z=\sqrt{r}\in \mathbb R_{+}$ one obtains
 $$\K=-\frac{1}{\pi  r \log(r^2)}+\pi ^{-1}\sum _{j=0}^{\infty } \left(\frac{r^{2 j}}{\left(-r^{1+2 j}+1\right)^2}+\frac{r^{2 j}}{\left(1-r^{2
j+1} \right)^2}\right) $$
$$\KX=\KY=\frac{1}{\pi  r^{\frac{3}{2}} \log(r^2)}+\pi ^{-1}\sum _{j=0}^{\infty } \left(\frac{2r^{\frac{1}{2}+4j}}{\left(-r^{1+2 j}+1\right)^3}-\frac{2r^{-\frac{1}{2}+2 j}}{\left(1-r^{2
j+1} \right)^3}\right) $$
$$\KXX=\KYY=-\frac{2}{\pi  r^2 \log(r^2)}+\pi ^{-1}\sum _{j=0}^{\infty } \left(
\frac{6r^{1+2j}}{\left(-r^{1+2 j}+1\right)^4}+
\frac{6r^{2j-1}}{\left(1-r^{2j+1}\right)^4}\right)$$
$$\KXY=-\frac{1}{\pi  r^2 \log(r^2)}+\pi ^{-1}\sum _{j=0}^{\infty } \Big(
\frac{6r^{1+6j}}{\left(-r^{1+2 j}+1\right)^4}+
\frac{2r^{4j}}{\left(-r^{1+2j}+1\right)^3}+$$
$$\frac{6r^{2j-1}}{\left(1-r^{1+2j}\right)^4}-
\frac{2r^{2j-1}}{\left(1-r^{2j+1}\right)^3}\Big)$$
$$\KXXY=\KXYY=\frac{2}{\pi  r^{\frac{5}{2}} \log(r^2)}+\pi ^{-1}\sum _{j=0}^{\infty }\Big( 
\frac{24r^{\frac{3}{2}+8j} }{\left(-r^{1+2 j}+1\right)^5}+
\frac{12r^{\frac{1}{2}+6j}}{\left(-r^{1+2j}+1\right)^4}-$$
$$
\frac{24r^{2j-\frac{3}{2}}}{\left(1-r^{1+2j}\right)^5}+
\frac{12r^{2j-\frac{3}{2}}}{\left(1-r^{2j+1}\right)^4}\Big)$$
$$\KXXYY=-\frac{4}{\pi  r^{3} \log(r^2)}+\pi ^{-1}\sum _{j=0}^{\infty } \Big(
\frac{120 r^{10j+2}}{\left(-r^{1+2 j}+1\right)^6}+
\frac{96 r^{1+8j}}{\left(-r^{1+2 j}+1\right)^5}+$$$$
\frac{12r^{6j}}{\left(-r^{1+2 j}+1\right)^4}\Big)
+\left(
\frac{120r^{2j-2}}{\left(1-r^{2j+1}\right)^6}-
\frac{96r^{2j-2}}{\left(1-r^{2j+1}\right)^5}+
\frac{12r^{2j-2}}{\left(1-r^{2j+1}\right)^4}\right)$$

All of the above series are locally uniformly convergent in the unit circle and the summands are of the form $f(\sqrt r)$, with $f$ real analytic (with exception for the very first summands, which may contribute some singular terms). Therefore each of the above expressions is of the form $F(\sqrt{r})+$ singular part, with $F$- real-analytic, hence it is justified to write:

$$\K=-\frac{1}{\pi  r \log(r^2)}+\pi ^{-1}\left(2+4r+8r^2+O(r^3)\right) $$
$$\KX=\KY=\frac{1}{\pi  r^{\frac{3}{2}} \log(r^2)}+\pi^{-1} \left(\frac{-2}{\sqrt{r}}-4\sqrt{r}-8r^{\frac{3}{2}}+O(r^{\frac{5}{2}})\right) $$
$$\KXX=\KYY=-\frac{2}{\pi  r^2 \log(r^2)}+\pi ^{-1}\left(\frac{6}{r}+24+O(r)\right)$$
$$\KXY=-\frac{1}{\pi  r^2 \log(r^2)}+\pi ^{-1}\left(\frac{4}{r}+20+64r+O(r^2)\right)$$
$$\KXXY=\KXYY=\frac{2}{\pi  r^{\frac{5}{2}} \log(r^2)}+\pi^{-1}\left(-\frac{12}{r^{\frac{3}{2}}}-\frac{72}{\sqrt{r}}+O(r^{\frac{1}{2}})\right)$$
$$\KXXYY=-\frac{4}{\pi  r^{3} \log(r^2)}+\pi ^{-1}\left(\frac{36}{r^{2}}+\frac{288}{r}+O(C)\right),$$
where as usual $O(r^{\alpha})$ is a substitute for an expression, which divided by $r^{\alpha}$ is bounded when $r$ tends to $0$, $r>0$. Note that in this representation always the powers of $r$ in the term $\frac{1}{\pi  r^{\alpha} \log(r^2)}$ and in $O(r^\beta)$ are in the following relation
\begin{align*}
 \tag{\dag}\alpha+\beta\geq3
\end{align*}

Our first task is to show that one can get rid of the $O$'s in the expressions representing $A_{j}(r)$, $j=1..24$ in the above notation. Let $A^{*}_{j}(r)$ be the expression $A_{j}(r)$, with the $O$'s deleted.

One has to show that $\lim_{r\to 0^{+}}A_{j}(r)-A^{*}_{j}(r)=0$, $j=1..24$

After some elementary algebraic manipulations one obtains
$$A_{j}(r)=\frac{\pm1}{(r\log(r^2))^{p_{j}}}\frac{\prod_{i=1}^{m_{j}}(a_{0}^{ij}+a_{1}^{ij}r^{\alpha_{1}^{ij}}\log(r^2)+..+a_{k(ij)}^{ij}r^{\alpha_{k(ij)}^{ij} }\log(r^2)O(r^{\gamma_{k(ij)}^{ij} }))}{\prod_{i=1}^{n_{j}}(b_{0}^{ij}+b_{1}^{ij}r^{\beta_{1}^{ij}}\log(r^2)+..+b_{s(ij)}^{ij}r^{\beta_{s(ij)}^{ij}}\log(r^2)O(r^{\delta_{s(ij)}^{ij} }))}$$

Here $p_{j}=2$ for $j=6,7,8,9,12,16,17,22,24$, otherwise $p_{j}=3$, the numbers $\alpha_{l}^{ij}$ and $\beta_{l}^{ij}$ form  ascending (with respect to $l$) sequences of $k(ij)$, respectively $s(ij)$  positive rational numbers.  The notation $k(ij)$ and $s(ij)$ is to stress the fact, that the length of each sequence depends on both $i$ and $j$. The $a_{l}^{ij}$ and $b_{l}^{ij}$ are some (rational) constants.  Finally $\alpha_{k(ij)}^{ij}+\gamma_{k(ij)}^{ij}\geq 3$ and $\beta_{s(ij)}^{ij}+\delta_{s(ij)}^{ij}\geq 3$ (by \dag, the reason for which $\alpha+\beta=4$, not $3$ it the expansions of $\K$,$\KX$,$\KXY$ is because the lowest power of $r$ disappears when one manipulates with $\mianown$). In what follows we  shorten $k(ij)$ and $s(ij)$ to just $k$ and $s$ but still keep in mind the mentioned dependence.

The numerator of $A_{j}(r)-A_{j}^{*}(r)$ is
$$\Big(\prod_{i=1}^{m_{j}}\Big(a_{0}^{ij}+
\sum_{l=1}^{k-1}a_{l}^{ij}r^{\alpha_{l}^{ij}}\log(r^2)+a_{k}^{ij}r^{\alpha_{k}^{ij} }\log(r^2)O(r^{\gamma_{k}^{ij}})\Big)\Big)\Big(\prod_{i=1}^{n_{j}}\Big(b_{0}^{ij}+
\sum_{l=1}^{s-1}(b_{l}^{ij}r^{\beta_{l}^{ij}}\log(r^2)\Big)\Big)$$$$-\Big(\prod_{i=1}^{n_{j}}\Big(b_{0}^{ij}+
\sum_{l=1}^{s-1}b_{l}^{ij}r^{\beta_{l}^{ij}}\log(r^2)+b_{s}^{ij}r^{\beta_{s}^{ij} }\log(r^2)O(r^{\delta_{s}^{ij}})\Big)\Big)\Big(\prod_{i=1}^{m_{j}}\Big(a_{0}^{ij}+
\sum_{l=1}^{k-1}(a_{l}^{ij}r^{\alpha_{l}^{ij}}\log(r^2)\Big)\Big)$$ and we see that all the terms not containing an $O$ kill each other. What remains is $O(r^3\log(r))$ (by \dag)

The denominator is

$$(r\log(r^2))^{p_{j}}\Big(\prod_{i=1}^{n_{j}}\Big(b_{0}^{ij}+
\sum_{l=1}^{s-1}b_{l}^{ij}r^{\beta_{l}^{ij}}\log(r^2)+b_{s}^{ij}r^{\beta_{s}^{ij} }\log(r^2)O(r^{\delta_{s}^{ij}})\Big)\Big)\times$$$$\Big(\prod_{i=1}^{n_{j}}\Big(b_{0}^{ij}+
\sum_{l=1}^{s-1}b_{l}^{ij}r^{\beta_{l}^{ij}}\log(r^2)\Big)\Big),$$

which is $O((r\log(r^2))^{p_{j}})$ (no $b_{0}^{ij}$ is zero), and hence the ratio tends to $0$.

The terms $A^{*}_{j}$ can be asymptotically evaluated as follows:
$$A^{*}_{23}(r)\approx-\frac{1}{2 r^3 \left(\log(r^2)\right)^3} +\frac{7}{r^2 \left(\log(r^2)\right)^2}+\frac{12}{ r^2 \left(\log(r^2)\right)^3}-\frac{30}{r
\left(\log(r^2)\right)}$$

$$  A^{*}_{21}(r)\approx \frac{1}{2 r^3 \left(\log(r^2)\right)^3}-\frac{6}{r^2
\left(\log(r^2)\right)^2}-\frac{12}{ r^2 \left(\log(r^2)\right)^3}+\frac{22}{r \left(\log(r^2)\right)}$$

$$  A^{*}_{20}(r)\approx \frac{1}{2 r^3 \left(\log(r^2)\right)^3}-\frac{11}{2 r^2 \left(\log(r^2)\right)^2}-\frac{12}{ r^2 \left(\log(r^2)\right)^3}+\frac{18}{r \left(\log(r^2)\right)}$$

$$ A^{*}_{19}(r)\approx-\frac{1}{2 r^3 \left(\log(r^2)\right)^3}+\frac{5}{r^2 \left(\log(r^2)\right)^2}+\frac{12}{ r^2 \left(\log(r^2)\right)^3}-\frac{14}{r
\left(\log(r^2)\right)}$$

$$  A^{*}_{18}(r)\approx \frac{1}{2 r^3 \left(\log(r^2)\right)^3}-\frac{11}{2
r^2 \left(\log(r^2)\right)^2}-\frac{12}{ r^2 \left(\log(r^2)\right)^3}+\frac{18}{r \left(\log(r^2)\right)}$$

$$  A^{*}_{15}(r)\approx-\frac{1}{2 r^3 \left(\log(r^2)\right)^3} +\frac{9}{2 r^2 \left(\log(r^2)\right)^2}+\frac{12}{
r^2 \left(\log(r^2)\right)^3}-\frac{13}{r \left(\log(r^2)\right)}$$

$$ A^{*}_{14}(r)\approx-\frac{1}{2 r^3 \left(\log(r^2)\right)^3} +\frac{4}{r^2 \left(\log(r^2)\right)^2}+\frac{12}{
r^2 \left(\log(r^2)\right)^3}-\frac{21}{2
r \left(\log(r^2)\right)}$$

$$  A^{*}_{13}(r)\approx \frac{1}{2 r^3 \left(\log(r^2)\right)^3}-\frac{7}{2
r^2 \left(\log(r^2)\right)^2}-\frac{12}{
r^2 \left(\log(r^2)\right)^3}+\frac{8}{r \left(\log(r^2)\right)}$$

$$  A^{*}_{11}(r)\approx \frac{1}{2 r^3 \left(\log(r^2)\right)^3}-\frac{6}{r^2 \left(\log(r^2)\right)^2}-\frac{12}{
r^2 \left(\log(r^2)\right)^3}+\frac{22}{r \left(\log(r^2)\right)}$$

$$  A^{*}_{10}(r)\approx-\frac{1}{2 r^3 \left(\log(r^2)\right)^3}+\frac{5}{r^2 \left(\log(r^2)\right)^2}+\frac{12}{
r^2 \left(\log(r^2)\right)^3}-\frac{14}{r
\left(\log(r^2)\right)}$$

$$  A^{*}_{5}(r)\approx-\frac{1}{2 r^3 \left(\log(r^2)\right)^3}+\frac{5}{r^2
\left(\log(r^2)\right)^2}+\frac{12}{
r^2 \left(\log(r^2)\right)^3}-\frac{16}{r \left(\log(r^2)\right)}$$

$$  A^{*}_{4}(r)\approx-\frac{1}{2 r^3 \left(\log(r^2)\right)^3} +\frac{9}{2 r^2 \left(\log(r^2)\right)^2}+\frac{12}{
r^2 \left(\log(r^2)\right)^3}-\frac{13}{r \left(\log(r^2)\right)}$$

$$ A^{*}_{3}(r)\approx \frac{1}{r^3 \left(\log(r^2)\right)^3} -\frac{8}{r^2 \left(\log(r^2)\right)^2}-\frac{24}{r^2 \left(\log(r^2)\right)^3}+\frac{20}{r
\left(\log(r^2)\right)}$$

$$  A^{*}_{2}(r)\approx \frac{1}{2 r^3 \left(\log(r^2)\right)^3} -\frac{7}{2
r^2 \left(\log(r^2)\right)^2}-\frac{12}{ r^2 \left(\log(r^2)\right)^3}+\frac{8}{r \left(\log(r^2)\right)}$$

$$  A^{*}_{1}(r)\approx-\frac{1}{2 r^3 \left(\log(r^2)\right)^3} +\frac{3}{r^2 \left(\log(r^2)\right)^2}+\frac{12}{ r^2 \left(\log(r^2)\right)^3}-\frac{6}{r \left(\log(r^2)\right)}, $$
where each time 
$$A^{*}_{j}(r)\approx \frac{a_{j}}{ r^3 \left(\log(r^2)\right)^3} +\frac{b_{j}}{r^2 \left(\log(r^2)\right)^2}+\frac{c_{j}}{ r^2 \left(\log(r^2)\right)^3}+\frac{d_{j}}{r \left(\log(r^2)\right)}$$
should read as
$$\lim_{r\to0^{+}}A^{*}_{j}(r)r^3(\log(r^2))^3=a_{j},\lim_{r\to0^{+}}(A^{*}_{j}(r)-\frac{a_{j}}{r^3(\log(r^2))^3})(r^2(\log(r^2))^2)=b_{j},$$

$$\lim_{r\to0^{+}}\Big(A^{*}_{j}(r)-\frac{a_{j}}{r^3(\log(r^2))^3}-\frac{b_{j}}{r^2(\log(r^2))^2}\Big)r^2(\log(r^2))^3=c_{j},$$

$$\lim_{r\to0^{+}}\Big(A^{*}_{j}(r)-\frac{a_{j}}{r^2(\log(r^2))^2}-\frac{b_{j}}{r\log(r^2)}-\frac{c_{j}}{r^2(\log(r^2))^3}\Big)r\log(r^2)=d_{j}.$$

and
$$A^{*}_{24}(r)\approx\frac{1}{r^2(\log(r^2))^2}+\frac{11}{r\log(r^2)}$$
$$A^{*}_{22}(r)\approx-\frac{1}{r^2(\log(r^2))^2}-\frac{8}{r\log(r^2)}$$
$$A^{*}_{17}(r)\approx-\frac{1}{2r^2(\log(r^2))^2}-\frac{6}{r\log(r^2)}$$
$$A^{*}_{16}(r)\approx\frac{1}{2r^2(\log(r^2))^2}+\frac{5}{r\log(r^2)}$$
$$A^{*}_{12}(r)\approx-\frac{1}{2r^2(\log(r^2))^2}-\frac{8}{r\log(r^2)}$$
$$A^{*}_{9}(r)\approx-\frac{1}{2r^2(\log(r^2))^2}-\frac{4}{r\log(r^2)}$$
$$A^{*}_{8}(r)\approx\frac{2}{r^2(\log(r^2))^2} +\frac{12}{r\log(r^2)}$$
$$A^{*}_{7}(r)\approx\frac{1}{2r^2(\log(r^2))^2} +\frac{5}{r\log(r^2)}$$
$$A^{*}_{6}(r)\approx-\frac{3}{2r^2(\log(r^2))^2} -\frac{6}{r\log(r^2)},$$
where each time 
$$A^{*}_{j}(r)\approx\frac{a_{j}}{r^2(\log(r^2))^2}+\frac{b_{j}}{r(\log(r^2))}$$
should read as
$$\lim_{r\to0^{+}}A^{*}_{j}(r)r^2(\log(r^2))^2=a_{j},\lim_{r\to0^{+}}(A^{*}_{j}(r)-\frac{a_{j}}{r^2(\log(r^2))^2})(r\log(r^2))=b_{j}.$$
Hence
$$\sum_{j=1}^{24}A^{*}_{j}(r)\approx \frac{0}{ r^3 \log (r^2)^3} +\frac{0}{r^2 \log (r^2)^2}+\frac{0}{ r^2 \log (r^2)^3}+\frac{1}{2r \log (r^2)}$$
$$\lim_{r\to 0^{+}}\sum_{j=1}^{24}A_{j}(r)=\lim_{r\to 0^{+}}\sum_{j=1}^{24}A^{*}_{j}(r)=-\infty$$

This proves the first part of Theorem \ref{glownetw}.

For the special choice $z=r^{\frac{3}{10}}\in \mathbb R_{+}$ one has
 $$\K=-\frac{1}{\pi  r^{\frac{3}{5}} \log(r^{2})}+\pi ^{-1}\sum _{j=0}^{\infty } \left(\frac{r^{2 j} }{\left(-r^{\frac{3}{5}+2 j}+1\right)^2}+\frac{r^{2 j+\frac{4}{5}} }{\left(1-r^{2
j+\frac{7}{5}} \right)^2}\right) $$
$$\KX=\frac{1}{\pi  r^{\frac{9}{10}} \log(r^2)}+\pi ^{-1}\sum _{j=0}^{\infty } \left(\frac{2r^{\frac{3}{10}+4j}}{\left(-r^{\frac{3}{5}+2 j}+1\right)^3}-\frac{2r^{\frac{1}{2}+2 j}}{\left(1-r^{2
j+\frac{7}{5}} \right)^3}\right) $$
$$\KY=\KX$$
$$\KXX=-\frac{2}{\pi  r^{\frac{6}{5}} \log(r^2)}+\pi ^{-1}\sum _{j=0}^{\infty } \left(
\frac{6r^{\frac{3}{5}+6j}}{\left(-r^{\frac{3}{5}+2 j}+1\right)^4}+
\frac{6r^{\frac{1}{5}+2j}}{\left(1-r^{2j+\frac{7}{5}}\right)^4}\right)$$
$$\KYY=\KXX$$
$$\KXY=-\frac{1}{\pi  r^{\frac{6}{5}} \log(r^2)}+\pi ^{-1}\sum _{j=0}^{\infty }\Big[ 
\frac{6r^{\frac{3}{5}+6j}}{\left(-r^{\frac{3}{5}+2 j}+1\right)^4}+
\frac{2r^{4j}}{\left(-r^{\frac{3}{5}+2j}+1\right)^3}+$$
$$\frac{6r^{2j+\frac{1}{5}}}{\left(1-r^{\frac{7}{5}+2j}\right)^4}-
\frac{2r^{2j+\frac{1}{5}}}{\left(1-r^{2j+\frac{7}{5}}\right)^3}\Big]$$

$$\KXXY=\frac{2}{\pi  r^{\frac{3}{2}} \log(r^2)}+\pi ^{-1}\sum _{j=0}^{\infty }\Big[ 
\frac{24r^{\frac{9}{10}+8j} }{\left(-r^{\frac{3}{5}+2 j}+1\right)^5}+
\frac{12r^{\frac{3}{10}+6j}}{\left(-r^{\frac{3}{5}+2j}+1\right)^4}-$$
$$
\frac{24r^{2j-\frac{1}{10}}}{\left(1-r^{\frac{7}{5}+2j}\right)^5}+
\frac{12r^{2j-\frac{1}{10}}}{\left(1-r^{2j+\frac{7}{5}}\right)^4}\Big]$$
$$\KXYY=\KXXY$$

$$\KXXYY=-\frac{4}{\pi  r^{\frac{9}{5}} \log(r^2)}+\pi ^{-1}\sum _{j=0}^{\infty } \Big[
\frac{120 r^{10j+\frac{6}{5}}}{\left(-r^{\frac{3}{5}+2 j}+1\right)^6}+
\frac{96 r^{\frac{3}{5}+8j}}{\left(-r^{\frac{3}{5}+2 j}+1\right)^5}+$$$$
\frac{12r^{6j}}{\left(-r^{\frac{3}{5}+2 j}+1\right)^4}
+\frac{120r^{2j-\frac{2}{5}}}{\left(1-r^{2j+\frac{7}{5}}\right)^6}-
\frac{96r^{2j-\frac{2}{5}}}{\left(1-r^{2j+\frac{7}{5}}\right)^5}+
\frac{12r^{2j-\frac{2}{5}}}{\left(1-r^{2j+\frac{7}{5}}\right)^4}\Big]$$

 As above each sum is of the form $G(r^\frac{1}{10})+$ singular part, with $G$- real-analytic. Now
$$\K=-\frac{1}{\pi  r^{\frac{3}{5}}\log(r^2)}+\pi ^{-1}(1+2r^{\frac{3}{5}}+ r^{\frac{4}{5}}+3r^{\frac{6}{5}}+O(r^{\frac{9}{5}})) $$
$$\KX=\KY=\frac{1}{\pi  r^{\frac{9}{10}} \log(r^2)}+\pi^{-1}(2r^{\frac{3}{10}}-2r^{\frac{1}{2}}+6r^{\frac{9}{10}}+O(r^{\frac{3}{2}})) $$
$$\KXX=\KYY=-\frac{2}{\pi  r^{\frac{6}{5}} \log(r^2)}+\pi^{-1}(6r^{\frac{1}{5}}+O(r^{\frac{3}{5}}))$$
$$\KXY=-\frac{1}{\pi  r^{\frac{6}{5}} \log(r^2)}+\pi ^{-1}(2+4r^{\frac{1}{5}}+12r^{\frac{3}{5}}+O(r{^\frac{6}{5}}))$$
$$\KXXY=\KXYY=\frac{2}{\pi  r^{\frac{3}{2}} \log(r^2)}+\pi^{-1}\left(-\frac{12}{r^{\frac{1}{10}}}+O(r^{\frac{3}{10}})\right)$$
$$\KXXYY=-\frac{4}{\pi  r^{\frac{9}{5}} \log(r^2)}+\pi ^{-1}\left(\frac{36}{r^{\frac{2}{5}}}+O(C)\right)$$

This time the goal is to have  
 \begin{align*}
\tag{\dag\dag}\alpha+\beta\geq\frac{9}{5},
\end{align*}
 for $\alpha$ and $\beta$ being
 the powers of $r$ in  $\frac{1}{\pi  r^{\alpha} \log(r^2)}$ and  $O(r^\beta)$ respectively. The argument with passing form $A_{j}(r)$ to $A^{*}_{j}(r)$ is almost the same, one only has to adjust the terms $(r\log(r^2))^{p_{j}}$ with $r^{p_{j}}\log(r^2)^{\frac{5}{3}p_{j}} $ and consider new $p_{j}$'s, namely $p_{j}=\frac{6}{5}$, for $j=6,7,8,9,12,16,17,22,24$ and $p_{j}=\frac{9}{5}$ otherwise.

Now the asymptotic expansions of $A_{j}^{*}(r)$ are:

$$A^{*}_{23}(r)\approx-\frac{4}{r^{\frac{9}{5}}\left(\log(r^2)\right)^3} +\frac{16}{r^{\frac{6}{5}}\left(\log(r^2)\right)^2}+\frac{96}{r^{\frac{6}{5}}\left(\log(r^2)\right)^3}+\frac{12}{r \left(\log(r^2)\right)^3}-\frac{24}{r^{\frac{3}{5}}
\left(\log(r^2)\right)}  $$$$ -\frac{376}{r^{\frac{3}{5}}\left(\log(r^2)\right)^2}-\frac{1212}{r^{\frac{3}{5}}\left(\log(r^2)\right)^3}-\frac{32}{r^{\frac{2}{5}}\left(\log(r^2)\right)^2}-\frac{384}{r^{\frac{2}{5}}\left(\log(r^2)\right)^3}-\frac{24}{r^{\frac{1}{5}}\left(\log(r^2)\right)^3} +16$$

$$ A^{*}_{21}(r)\approx\frac{4}{r^{\frac{9}{5}}\left(\log(r^2)\right)^3}-\frac{12}{r^{\frac{6}{5}}\left(\log(r^2)\right)^2}-\frac{96}{r^{\frac{6}{5}}\left(\log(r^2)\right)^3}-\frac{12}{r
\left(\log(r^2)\right)^3}+\frac{12}{r^{\frac{3}{5}}\left(\log(r^2)\right)} $$$$  +\frac{288}{r^{\frac{3}{5}}\left(\log(r^2)\right)^2}+\frac{1212}{r^{\frac{3}{5}}\left(\log(r^2)\right)^3} +\frac{24}{r^{\frac{2}{5}}
\left(\log(r^2)\right)^2} +\frac{384}{r^{\frac{2}{5}}\left(\log(r^2)\right)^3}+\frac{24}{r^{\frac{1}{5}}
\left(\log(r^2)\right)^3}-4$$

$$ A^{*}_{20}(r)\approx\frac{4}{r^{\frac{9}{5}}\left(\log(r^2)\right)^3}-\frac{12}{r^{\frac{6}{5}}\left(\log(r^2)\right)^2}-\frac{96}{r^{\frac{6}{5}}
\left(\log(r^2)\right)^3}-\frac{12}{r \left(\log(r^2)\right)^3}+\frac{12}{r^{\frac{3}{5}}\left(\log(r^2)\right)}  $$$$  +\frac{296}{r^{\frac{3}{5}}\left(\log(r^2)\right)^2}+\frac{1212}{r^{\frac{3}{5}}\left(\log(r^2)\right)^3}+\frac{28}{r^{\frac{2}{5}
}\left(\log(r^2)\right)^2}+\frac{384}{r^{\frac{2}{5}}\left(\log(r^2)\right)^3} +\frac{24}{r^{\frac{1}{5}
}\left(\log(r^2)\right)^3}-4$$

$$ A^{*}_{19}(r)\approx-\frac{4}{r^{\frac{9}{5}}\left(\log(r^2)\right)^3}+\frac{8}{r^{\frac{6}{5}}\left(\log(r^2)\right)^2}+\frac{96}{r^{\frac{6}{5}
}\left(\log(r^2)\right)^3}+\frac{12}{r \left(\log(r^2)\right)^3}-\frac{4}{r^{\frac{3}{5}
}\left(\log(r^2)\right)}  $$$$ -\frac{224}{r^{\frac{3}{5}}\left(\log(r^2)\right)^2}-\frac{1212}{r^{\frac{3}{5}}\left(\log(r^2)\right)^3}-\frac{16}{r^{\frac{2}{5}}\left(\log(r^2)\right)^2}-\frac{384}{r^{\frac{2}{5}}\left(\log(r^2)\right)^3} -\frac{24}{r^{\frac{1}{5}
}\left(\log(r^2)\right)^3}$$

$$ A^{*}_{18}(r)\approx\frac{4}{r^{\frac{9}{5}}\left(\log(r^2)\right)^3}-\frac{12}{r^{\frac{6}{5}}\left(\log(r^2)\right)^2}-\frac{96}{r^{\frac{6}{5}}\left(\log(r^2)\right)^3}-\frac{12}{r
\left(\log(r^2)\right)^3}+\frac{12}{r^{\frac{3}{5}}\left(\log(r^2)\right)}  $$$$  +\frac{296}{r^{\frac{3}{5}}\left(\log(r^2)\right)^2}+\frac{1212}{r^{\frac{3}{5}}\left(\log(r^2)\right)^3}+\frac{28}{r^{\frac{2}{5}
}\left(\log(r^2)\right)^2}+\frac{384}{r^{\frac{2}{5}}\left(\log(r^2)\right)^3} +\frac{24}{r^{\frac{1}{5}
}\left(\log(r^2)\right)^3}-4$$

$$ A^{*}_{15}(r)\approx-\frac{4}{r^{\frac{9}{5}}\left(\log(r^2)\right)^3}+\frac{8}{r^{\frac{6}{5}}\left(\log(r^2)\right)^2}+\frac{96}{r^{\frac{6}{5}
}\left(\log(r^2)\right)^3}+\frac{12}{r \left(\log(r^2)\right)^3}-\frac{4}{r^{\frac{3}{5}
}\left(\log(r^2)\right)}  $$$$  -\frac{208}{r^{\frac{3}{5}}\left(\log(r^2)\right)^2}-\frac{1212}{r^{\frac{3}{5}}\left(\log(r^2)\right)^3}-\frac{20}{r^{\frac{2}{5}}\left(\log(r^2)\right)^2}-\frac{384}{r^{\frac{2}{5}}\left(\log(r^2)\right)^3} -\frac{24}{r^{\frac{1}{5}
}\left(\log(r^2)\right)^3}$$

$$ A^{*}_{14}(r)\approx-\frac{4}{r^{\frac{9}{5}}\left(\log(r^2)\right)^3}+\frac{8}{r^{\frac{6}{5}}\left(\log(r^2)\right)^2}+\frac{96}{r^{\frac{6}{5}}\left(\log(r^2)\right)^3}+\frac{12}{r
\left(\log(r^2)\right)^3}-\frac{4}{r^{\frac{3}{5}
}\left(\log(r^2)\right)}  $$$$  -\frac{216}{r^{\frac{3}{5}}\left(\log(r^2)\right)^2}-\frac{1212}{r^{\frac{3}{5}}\left(\log(r^2)\right)^3}-\frac{24}{r^{\frac{2}{5}}\left(\log(r^2)\right)^2}-\frac{384}{r^{\frac{2}{5}}\left(\log(r^2)\right)^3} -\frac{24}{r^{\frac{1}{5}
}\left(\log(r^2)\right)^3}$$

$$ A^{*}_{13}(r)\approx\frac{4}{r^{\frac{9}{5}}\left(\log(r^2)\right)^3}-\frac{4}{r^{\frac{6}{5}}\left(\log(r^2)\right)^2}-\frac{96}{r^{\frac{6}{5}}\left(\log(r^2)\right)^3}-\frac{12}{r
\left(\log(r^2)\right)^3}  $$$$  +\frac{144}{r^{\frac{3}{5}}\left(\log(r^2)\right)^2}+\frac{1212}{r^{\frac{3}{5}}\left(\log(r^2)\right)^3}+\frac{12}{r^{\frac{2}{5}}\left(\log(r^2)\right)^2}+\frac{384}{r^{\frac{2}{5}}\left(\log(r^2)\right)^3}+\frac{24}{r^{\frac{1}{5}
}\left(\log(r^2)\right)^3}$$

$$ A^{*}_{11}(r)\approx\frac{4}{r^{\frac{9}{5}}\left(\log(r^2)\right)^3}-\frac{12}{r^{\frac{6}{5}
}\left(\log(r^2)\right)^2}-\frac{96}{r^{\frac{6}{5}}\left(\log(r^2)\right)^3}-\frac{12}{r \left(\log(r^2)\right)^3}+\frac{12}{r^{\frac{3}{5}
}\left(\log(r^2)\right)}  $$$$  +\frac{288}{r^{\frac{3}{5}}\left(\log(r^2)\right)^2}+\frac{1212}{r^{\frac{3}{5}
}\left(\log(r^2)\right)^3}+\frac{24}{r^{\frac{2}{5}}\left(\log(r^2)\right)^2} +\frac{384}{r^{\frac{2}{5}}\left(\log(r^2)\right)^3}+\frac{24}{r^{\frac{1}{5}}\left(\log(r^2)\right)^3}-4$$

$$ A^{*}_{10}(r)\approx-\frac{4}{r^{\frac{9}{5}}\left(\log(r^2)\right)^3}+\frac{8}{r^{\frac{6}{5}}\left(\log(r^2)\right)^2}+\frac{96}{r^{\frac{6}{5}}\left(\log(r^2)\right)^3}+\frac{12}{r
\left(\log(r^2)\right)^3}-\frac{4}{r^{\frac{3}{5}
}\left(\log(r^2)\right)}  $$$$  -\frac{224}{r^{\frac{3}{5}}\left(\log(r^2)\right)^2}-\frac{1212}{r^{\frac{3}{5}}\left(\log(r^2)\right)^3}-\frac{16}{r^{\frac{2}{5}}\left(\log(r^2)\right)^2}-\frac{384}{r^{\frac{2}{5}}\left(\log(r^2)\right)^3} -\frac{24}{r^{\frac{1}{5}
}\left(\log(r^2)\right)^3}$$

$$ A^{*}_{5}(r)\approx-\frac{4}{r^{\frac{9}{5}}\left(\log(r^2)\right)^3}+\frac{8}{r^{\frac{6}{5}}\left(\log(r^2)\right)^2}+\frac{96}{r^{\frac{6}{5}}\left(\log(r^2)\right)^3}+\frac{12}{r
\left(\log(r^2)\right)^3}-\frac{4}{r^{\frac{3}{5}}
\left(\log(r^2)\right)} $$$$  -\frac{200}{r^{\frac{3}{5}}\left(\log(r^2)\right)^2}-\frac{1212}{r^{\frac{3}{5}}\left(\log(r^2)\right)^3}-\frac{16}{r^{\frac{2}{5}}\left(\log(r^2)\right)^2} -\frac{384}{r^{\frac{2}{5}}\left(\log(r^2)\right)^3}-\frac{24}{r^{\frac{1}{5
}}\left(\log(r^2)\right)^3}$$

$$ A^{*}_{4}(r)\approx-\frac{4}{r^{\frac{9}{5}}\left(\log(r^2)\right)^3}+\frac{8}{r^{\frac{6}{5}}\left(\log(r^2)\right)^2}+\frac{96}{r^{\frac{6}{5}}\left(\log(r^2)\right)^3}+\frac{12}{r
\left(\log(r^2)\right)^3}-\frac{4}{r^{\frac{3}{5
}}\left(\log(r^2)\right)} $$$$ -\frac{208}{r^{\frac{3}{5}}
\left(\log(r^2)\right)^2}-\frac{1212}{r^{\frac{3}{5}}
\left(\log(r^2)\right)^3}-\frac{20}{r^{\frac{2}{5}}\left(\log(r^2)\right)^2}  -\frac{384}{r^{\frac{2}{5}
}\left(\log(r^2)\right)^3}-\frac{24}{r^{\frac{1}{5}
}\left(\log(r^2)\right)^3}$$

$$ A^{*}_{3}(r)\approx\frac{8}{r^{\frac{9}{5}}\left(\log(r^2)\right)^3} -\frac{8}{r^{\frac{6}{5}}\left(\log(r^2)\right)^2}-\frac{192}{r^{\frac{6}{5}}\left(\log(r^2)\right)^3}-\frac{24}{r
\left(\log(r^2)\right)^3}  $$$$ +\frac{272}{r^{\frac{3}{5}}\left(\log(r^2)\right)^2}+\frac{2424}{r^{\frac{3}{5}}\left(\log(r^2)\right)^3}+\frac{16}{r^{\frac{2}{5}}\left(\log(r^2)\right)^2}+\frac{768}{r^{\frac{2}{5}}\left(\log(r^2)\right)^3}+\frac{48}{r^{\frac{1}{5}
}\left(\log(r^2)\right)^3}$$

$$ A^{*}_{2}(r)\approx\frac{4}{r^{\frac{9}{5}}\left(\log(r^2)\right)^3}-\frac{4}{r^{\frac{6}{5}}\left(\log(r^2)\right)^2}-\frac{96}{r^{\frac{6}{5}}\left(\log(r^2)\right)^3}-\frac{12}{r \left(\log(r^2)\right)^3}  $$$$ +\frac{144}{r^{\frac{3}{5
}}\left(\log(r^2)\right)^2}+\frac{1212}{r^{\frac{3}{5
}}\left(\log(r^2)\right)^3}+\frac{12}{r^{\frac{2}{5}}\left(\log(r^2)\right)^2}+\frac{384}{r^{\frac{2}{5}}\left(\log(r^2)\right)^3}+\frac{24}{r^{\frac{1}{5}}\left(\log(r^2)\right)^3}$$

$$ A^{*}_{1}(r)\approx-\frac{4}{r^{\frac{9}{5}}\left(\log(r^2)\right)^3}+\frac{96}{r^{\frac{6}{5}
}\left(\log(r^2)\right)^3}+\frac{12}{r \left(\log(r^2)\right)^3}  $$$$ -\frac{72}{r^{\frac{3}{5}}\left(\log(r^2)\right)^2}-\frac{1212}{r^{\frac{3}{5}}\left(\log(r^2)\right)^3}-\frac{384}{r^{\frac{2}{5}}\left(\log(r^2)\right)^3}-\frac{24}{r^{\frac{1}{5
}}\left(\log(r^2)\right)^3}.$$

Again
$$A^{*}_{j}(r)\approx\frac{a_{j}}{r^{\frac{9}{5}}\left(\log(r^2)\right)^3} +\frac{b_{j}}{r^{\frac{6}{5}}\left(\log(r^2)\right)^2}+\frac{c_{j}}{r^{\frac{6}{5}}\left(\log(r^2)\right)^3}+\frac{d_{j}}{r \left(\log(r^2)\right)^3}+\frac{e_{j}}{r^{\frac{3}{5}}
\left(\log(r^2)\right)}  $$$$ +\frac{f_{j}}{r^{\frac{3}{5}}\left(\log(r^2)\right)^2}+\frac{g_{j}}{r^{\frac{3}{5}}\left(\log(r^2)\right)^3}+\frac{h_{j}}{r^{\frac{2}{5}}\left(\log(r^2)\right)^2}+\frac{i_{j}}{r^{\frac{2}{5}}\left(\log(r^2)\right)^3}+\frac{k_{j}}{r^{\frac{1}{5}}\left(\log(r^2)\right)^3} +l_{j}$$
means that
$$\lim_{r\to 0^{+}}A^{*}_{j}(r)r^{\frac{9}{5}}\left(\log(r^2)\right)^3=a_{j}, \ \  \lim_{r\to 0^{+}}(A^{*}_{j}(r)-\frac{a_{j}}{r^{\frac{9}{5}}\left(\log(r^2)\right)^3})r^{\frac{6}{5}}\left(\log(r^2)\right)^2=b_{j}$$
....
$$\lim_{r\to 0^{+}}\Big( A^{*}_{j}(r)-\frac{a_{j}}{r^{\frac{9}{5}}\left(\log(r^2)\right)^3} -\frac{b_{j}}{r^{\frac{6}{5}}\left(\log(r^2)\right)^2}-\frac{c_{j}}{r^{\frac{6}{5}}\left(\log(r^2)\right)^3}-\frac{d_{j}}{r \left(\log(r^2)\right)^3}-\frac{e_{j}}{r^{\frac{3}{5}}
\left(\log(r^2)\right)}  $$$$ -\frac{f_{j}}{r^{\frac{3}{5}}\left(\log(r^2)\right)^2}-\frac{g_{j}}{r^{\frac{3}{5}}\left(\log(r^2)\right)^3}-\frac{h_{j}}{r^{\frac{2}{5}}\left(\log(r^2)\right)^2}-\frac{i_{j}}{r^{\frac{2}{5}}\left(\log(r^2)\right)^3}-\frac{k_{j}}{r^{\frac{1}{5}}\left(\log(r^2)\right)^3}\Big)= l_{j}$$
(to simplify the calculations one can put $r=q^{10}$ and then find $\lim_{q\to 0} ...$)

The expansions of the other terms are:

$$A^{*}_{24}(r)\approx-\frac{4}{r^{\frac{6}{5}} \left(\log(r^2)\right)^2}+\frac{12}{r^{\frac{3}{5}}
\left(\log(r^2)\right)}+\frac{64}{r^{\frac{3}{5}} \left(\log(r^2)\right)^2}+\frac{8}{r^{\frac{2}{5}} \left(\log(r^2)\right)^2} -12$$

$$A^{*}_{22}(r)\approx\frac{4}{r^{\frac{6}{5}} \left(\log(r^2)\right)^2}-\frac{8}{r^{\frac{3}{5}} \left(\log(r^2)\right)}-\frac{64}{r^{\frac{3}{5}} \left(\log(r^2)\right)^2}-\frac{8}{r^{\frac{2}{5}}
\left(\log(r^2)\right)^2} +4$$

$$A^{*}_{17}(r)\approx\frac{4}{r^{\frac{6}{5}} \left(\log(r^2)\right)^2}-\frac{8}{r^{\frac{3}{5}} \left(\log(r^2)\right)}-\frac{64}{r^{\frac{3}{5}}
\left(\log(r^2)\right)^2}-\frac{8}{r^{\frac{2}{5}} \left(\log(r^2)\right)^2} +4$$

$$A^{*}_{16}(r)\approx-\frac{4}{r^{\frac{6}{5}} \left(\log(r^2)\right)^2}+\frac{4}{r^{\frac{3}{5}}
\left(\log(r^2)\right)}+\frac{64}{r^{\frac{3}{5}} \left(\log(r^2)\right)^2}+\frac{8}{r^{\frac{2}{5}} \left(\log(r^2)\right)^2} $$

$$A^{*}_{12}(r)\approx\frac{4}{r^{\frac{6}{5}} \left(\log(r^2)\right)^2}-\frac{8}{r^{\frac{3}{5}} \left(\log(r^2)\right)}-\frac{64}{r^{\frac{3}{5}} \left(\log(r^2)\right)^2}-\frac{8}{r^{\frac{2}{5}}
\left(\log(r^2)\right)^2}+4 $$

$$A^{*}_{9}(r)\approx\frac{2}{r^{\frac{6}{5}} \left(\log(r^2)\right)^2}-\frac{4}{r^{\frac{3}{5}} \left(\log(r^2)\right)}-\frac{32}{r^{\frac{3}{5}}
\left(\log(r^2)\right)^2}-\frac{4}{r^{\frac{2}{5}} \left(\log(r^2)\right)^2}+2 $$

$$A^{*}_{8}(r)\approx-\frac{8}{r^{\frac{6}{5}} \left(\log(r^2)\right)^2}+\frac{8}{r^{\frac{3}{5}}
\left(\log(r^2)\right)}+\frac{128}{r^{\frac{3}{5}} \left(\log(r^2)\right)^2}+\frac{16}{r^{\frac{2}{5}} \left(\log(r^2)\right)^2} $$

$$A^{*}_{7}(r)\approx-\frac{4}{r^{\frac{6}{5}} \left(\log(r^2)\right)^2}+\frac{4}{r^{\frac{3}{5}} \left(\log(r^2)\right)}+\frac{64}{r^{\frac{3}{5}} \left(\log(r^2)\right)^2}+\frac{8}{r^{\frac{2}{5}}
\left(\log(r^2)\right)^2} $$

$$A^{*}_{6}(r)\approx\frac{6}{r^{\frac{6}{5}} \left(\log(r^2)\right)^2}-\frac{96}{r^{\frac{3}{5}}
\left(\log(r^2)\right)^2} -\frac{12}{r^{\frac{2}{5}} \left(\log(r^2)\right)^2}. $$

Adding up one receives

$$\sum_{j=1}^{24} A_{j}^{*}(r)\approx 2.$$

Therefore
$$\lim_{r\to 0^{+}}A_{j}(r)=2.$$
\begin{observation}One has that
 $$R_{P_{r}}(e^{i\theta}z)=R_{P_{r}}(z),$$
for all $\theta\in \text{(}0,2\pi\text{]}$ and 
$$R_{P_{r}}(z)=R_{P_{r}}\left(\frac{r}{z}\right),$$
\end{observation}
One easily checks that both $z\to e^{i\theta}z$ and $z\to\frac{r}{z}$ are holomorphic automorphisms of $P_{r}$. Now everything follows from the invariance property of the holomorphic sectional curvature of the Bergman metric. One moreover sees, that for the choice $z=r^{\frac{7}{10}}=\frac{r}{r^{\frac{3}{10}}}$ the equality
$$\lim_{r\to 0^{+}} R_{P_{r}}(r^{\frac{7}{10}})=2$$
is also true, something which the figures also depict.
\begin{observation}\label{biholom} One has that
 $$R_{\{\rho_{1}<|w-z_{0}|<\rho_{2}\}}(z)=R_{P_{\frac{\rho_{1}}{\rho_{2}}}}\left(\frac{z-z_{0}}{\rho_{2}}\right),$$
\end{observation}
This is also a simple consequence of the fact, that $z\to \frac{z-z_{0}}{\rho_{2}}$ is biholomorphic between the two domains.
\end{section}
\begin{section}{An example}

Let $\{R_{j}\}_{j=1}^{\infty},\{r_{j}\}_{j=1}^{\infty},\{s_{j}\}_{j=1}^{\infty}$ be three sequences of positive real numbers, that obey the following conditions
\begin{align*}
\tag{i} &\sum_{j=1}^{\infty}R_{j}<\infty\\
\tag{ii} &r_{1}<\frac{R_{1}}{2}, \frac{r_{j}}{R_{j}} \text{\ \ is  decreasing and\ } lim_{j\to \infty}\frac{r_{j}}{R_{j}}=0\\
\tag{iii}  &s_{j}< min\Big\{2R_{j}\sin(0.07\pi), 2R_{j+1}  \sin(0.07\pi), R_{j}-\Big(\frac{r_{j}}{R_{j}}\Big)^{\frac{3}{10}},\\ R_{j+1}-\Big(\frac{r_{j+1}}{R_{j+1}}\Big)^{\frac{3}{10}}\Big\}
\end{align*}
Consider the following domain:
$\Omega=\bigcup_{j=1}^{\infty} \Omega_{j}$, where $\Omega_{1}$ is the circular annulus $\{z:r_{1}<|z|<R_{1}\}$, $\Omega_{j}$ is also an annulus with inner radius $r_{j}$ and outer radius $R_{j}$, which is centered on the positive real axis, on the right of $\Omega_{j-1}$ and overlaps with $\Omega_{j-1}$ in such a way, that the segment joining the two intersection points of the circles with radii $R_{j}$ and $R_{j-1}$ has length $s_{j-1}$.

By $i)$, $\Omega$ is bounded.
 
 Let $\Omega'_{j}=\Omega_{j}\setminus(K_{1}\cup K_{2})$, where 
$K_{1}$ is a circle, centered at the midpoint of the segment joining the intersection points of the circles with radii $R_{j}$ and $R_{j-1}$ (the outer 
boundaries of the annuli $\Omega_{j}$ and $\Omega_{j-1}$) and the radius of $K_{1}$ is $\frac{s_{j-1}}{2}$. $K_{2}$  has radius equal to $\frac{s_{j}}{2}$ and is  centered at the midpoint of the segment joining the intersection points of the circles with radii $R_{j}$ and $R_{j+1}$.

\begin{figure}[tbph]
\centering
\includegraphics[width=4.31 in, ] {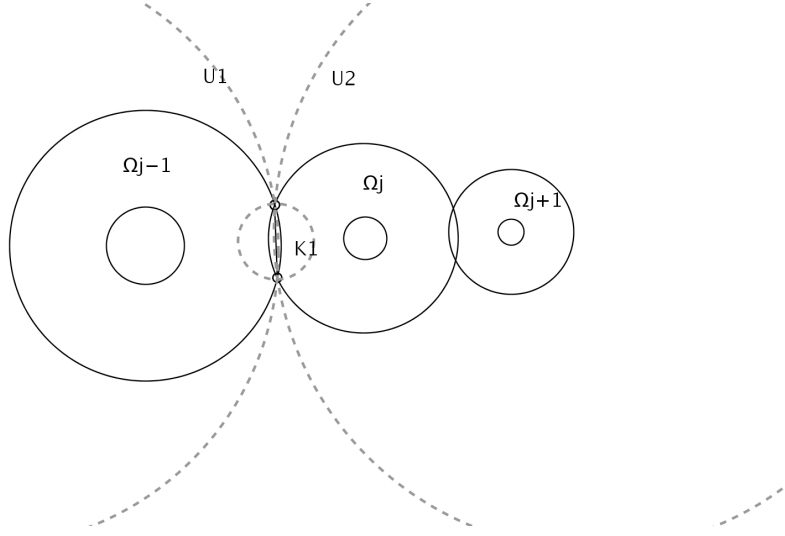}
\caption{Part of the domain $\Omega$}
\label{omega}
\end{figure}

 We begin with a lemma

\begin{lemma}\label{green} Let $z_{0}\in \Omega'_{j}$ then the sublevelset of the Green function $G_{\Omega}(z,z_{0})$,
$$\{G_{\Omega}(z,z_{0})<-1\}$$ 
is entirely contained in $\Omega_{j}$.
\end{lemma}
\begin{proof}
The Green function is obviously decreasing with respect to domain inclusions.Therefore it is enough to show that $\{G_{U}(z,z_{0})<-1\}\subset\Omega_{j}$, for some $\Omega\subset U$.

For the sake of simplicity one can translate $\Omega$, such that the upper intersection point of the circles with radii $R_{j-1}$ and $R_{j}$ is $0$. 
Choose  two circles $U_{1}$ and $U_{2}$, with radii $\rho_{1}$ and $\rho_{2}$, which intersect at $0$ and $-is_{j-1}$, such that $\Omega\subset U_{1}\cup U_{2}$. Clearly $\rho_{1}\geq R_{j-1}$ and $\rho_{2}>R_{j}$. The function $G_{U_{1}\cup U_{2}}(z,z_{0})$ can be explicitly calculated as $h\circ f$,where
$$h(w):=\log\left|\frac{w-f(z_{0})}{1-w\overline {f(z_{0})}}\right|$$
and 
$$f(z):=\frac{\left(\Big(\frac{1}{z}-\frac{i}{s_{j-1}}\Big)e^{-i\frac{\pi}{2}}e^{i\frac{\beta-\alpha}{2}}\right)^{\frac{\pi}{2\pi-\alpha-\beta}}-1}{\left(\Big(\frac{1}{z}-\frac{i}{s_{j-1}}\Big)e^{-i\frac{\pi}{2}}e^{i\frac{\beta-\alpha}{2}}\right)^{\frac{\pi}{2\pi-\alpha-\beta}}+1}$$
is the mapping that transforms $U_{1}\cup U_{2}$ conformally into the unit circle.Here $\alpha=\arcsin\frac{s_{j-1}}{2\rho_{2}}$, $\beta=\arcsin\frac{s_{j-1}}{2\rho_{1}}$.

The image of $U_{2}$ under $f$ is the intersection of the unit circle with a circle centered at the negative imaginary axis, passing through $\{1\}$ and $\{-1\}$ and such that the angle between it's tangent line at $\{-1\}$ and the line $-1+it$ is exactly $\frac{\pi}{2\pi-\beta-\alpha}\frac{\beta+\alpha}{2}$. The image of $U_{1}$ is exactly the conjugate of this set.
Now $f(\{G_{U_{1}\cup U_{2}}(z,z_{0})<1\})=\{w:\left|\frac{w-f(z_{0})}{1-w\overline {f(z_{0})}}\right|<e^{-1}\}$, which is the circle

$$\left|w-f(z_{0})\frac{1-e^{-2}}{1-e^{-2}|f(z_{0})|^2}\right|<\frac{e^{-1}(1-|f(z_{0})|^2)}{{1-e^{-2}|f(z_{0})|^2}}.$$

Now if one has that  $$z_{0}\in U_{2}, \text{\ \ \ }|Im f(z_{0})|\frac{1-e^{-2}}{1-e^{-2}|f(z_{0})|^2}\geq\frac{e^{-1}(1-|f(z_{0})|^2)}{{1-e^{-2}|f(z_{0})|^2}}$$ then this circle will stay in the lower halfdisc (and hence in $f(U_{2})$)

This inequality transforms easily into
\begin{equation}\label{good}\left|f(z_{0})-i\frac{1-e^{-2}}{2e^{-1}}\right|\geq\frac{1+e^{-2}}{2e^{-1}}, Im f(z_{0})<0.\end{equation}

So this is the set enclosed by the arcs of the unit circle and a circle that passes through $\{-1\}$ and $\{1\}$ and the angle between the real axis and the tangent line at $\{-1\}$ is $\arccos\frac{1-e^{-2}}{1+e^{-2}}\approx 0.22\pi$, $<0.23\pi$ (or minus this quantity if we refer to orientation).

On the other hand the image of the circle $|z+i\frac{s_{j-1}}{2}|<\frac{s_{j-1}}{2} (K_{1})$ is the set enclosed by two circular arcs joining $\{-1\}$ with $\{1\}$ characterized by the angle between the real axis ant the tangent line at $\{-1\}$ ($\frac{\pi-\beta+\alpha}{2}\frac{\pi}{2\pi-\beta-\alpha}$ and $(-)\frac{\pi+\beta-\alpha}{2}\frac{\pi}{2\pi-\beta-\alpha}$ respectively)
we see that if both $\alpha$ and $\beta$ are smaller than $0.07\pi$ then the slope of the lower arc is greater than $0.23\pi$ and hence lies in the region, defined by \ref{good}. Since $z_{0}\notin K_{1}$ then $f(z_{0})$ must lie below the discussed arc and hence in the desired region.

It remains to observe, that the condition $\alpha,\beta<0.07\pi$ is fulfilled by the choice of $s_{j-1}$ and that one can carry out the same argument for $K_{2}$.
\end{proof}
Let $\zeta$ be the rightmost boundary point of $\Omega$ ( the accumulation point of the annuli). Now one can give a proof of the main theorem 

\begin{proof}

It is clear that  $J_{i,\Omega}(z)\leq J_{i,\Omega_{j}}(z)$, for all $z\in \Omega_{j}, i=1,2,3$ and for all $j$ (see (\ref{j0}) ,(\ref{j1}), (\ref{j2}) ).

Let $z_{0}\in \Omega'_{j}$ and let $f_{i}(z)$ be the corresponding function, that realizes the supremum in the definition of $J_{i,\Omega_{j}}$. We have that 
$$f_{i}(z)\in \mathcal O(\Omega_{j})\cap L^2(\Omega_{j}),\int_{\Omega_{j}}|f_{i}|^2\leq 1, f_{i}^{(k)}(z_{0})=0, k=0,..,i-1.$$

Let $\chi$ be a real smooth function of a real variable, such that $\chi(x)=0$, for $x>-1$, $\chi(x)=1$, for $x<-2$, $0\leq\chi(x)\leq1$, for $-2\leq x\leq -1$, and $|\chi'(x)|<C$, globally for some  positive constant $C$.

By Lemma \ref{green} the $(0,1)$ form $\bar{\partial}( \chi\circ G_{\Omega}(z,z_{0})) . f_{i}(z)$ can be extended (trivially) to a smooth ($\bar\partial$-closed) form on the whole $\Omega$.

Note that $e^{G_{\Omega}(z,z_{0})}$ is a subharmonic function that satisfies $$\frac{\partial^2}{\partial z\partial \bar z}e^{G_{\Omega}(z,z_{0})}\geq\left|\frac{\partial}{\partial z}e^{G_{\Omega}(z,z_{0})}\right|^2,$$
in the weak sense. Therefore by the Donnelly-Fefferman estimate (see \cite{MR1415958} and especially  \cite{MR2139520}, where the passing from smooth to nonsmooth data is presented very clearly) one can find a solution $v$ of the $\bar\partial$ problem $$\bar\partial v_{i}=\bar{\partial}( \chi\circ G_{\Omega}(z,z_{0})) . f_{i}(z)$$ in $\Omega$, with

$$\int_{\Omega}|v_{i}|^2e^{-2(i+1)G_{\Omega}(z,z_{0})}\leq C'\int_{\Omega}\frac{|\bar{\partial}( \chi\circ G_{\Omega}(z,z_{0})) . f_{i}(z)|^2}{\frac{\partial^2}{\partial z\partial \bar z}e^{G_{\Omega}(z,z_{0})}}e^{-2(i+1)G_{\Omega}(z,z_{0})}\leq$$
$$C'\int_{\{-2<G_{\Omega}(z,z_{0})<-1\}}\frac{C^2|f_{i}|^2}{e^{(2i+3)G_{\Omega}(z,z_{0})}}\leq C'C^2 e^{2(2i+3)},$$
where $C'$ is a universal constant.

Moreover $v_{i}$ is holomorphic in a neighbourhood of $z_{0}$, and the above inequality ensures that $v_{i}^{(k)}(z_{0})=0, k=0,..,i$

The function $g_{i}=( \chi\circ G_{\Omega}(z,z_{0})) . f_{i}(z)-v_{i}$ is holomorphic in $\Omega$, agrees to the $i$-th order with $f_{i}$ at $z_{0}$ and
$$(\int_{\Omega}|g_{i}|^2)^{\frac{1}{2}}\leq (\int_{\Omega}|v_{i}|^2)^{\frac{1}{2}}+(\int_{\Omega}|( \chi\circ G_{\Omega}(z,z_{0})) . f_{i}(z)|^2)^{\frac{1}{2}}$$$$\leq (\int_{\Omega}|v_{i}|^2e^{-2(i+1)G_{\Omega}(z,z_{0})})^{\frac{1}{2}}+1\leq 1+\sqrt{C'C^2 e^{2(2i+3)}}$$

The choice of the function $\frac{g_{i}(z)}{1+\sqrt{C'C^2 e^{2(2i+3)}}}$ shows that 
$$J_{i,\Omega}(z)\geq \frac{J_{i,\Omega_{j}}(z)}{1+\sqrt{C'C^2 e^{2(2i+3)}}},$$ for all $z\in \Omega_{j}'$,
which is independent of $j$.

Hence
$$R_{\Omega}(z)=2-\frac{J_{0,\Omega}(z)J_{2,\Omega}(z)}{J_{1,\Omega}(z)^2}\leq 2-\frac{J_{0,\Omega_{j}}(z)J_{2,\Omega_{j}}(z)}{(1+\sqrt{C'C^2 e^{14}})(1+\sqrt{C'C^2 e^{6}})J_{1,\Omega_{j}}(z)^2}$$
$$=\frac{R_{\Omega_{j}}(z)+2C''-2}{C''}$$

and $$R_{\Omega}(z)\geq 2-\frac{(1+\sqrt{C'C^2 e^{10}})^2J_{0,\Omega_{j}}(z)J_{2,\Omega_{j}}(z)}{J_{1,\Omega_{j}}(z)^2}=C'''(R_{\Omega_{j}}(z)-2+\frac{2}{C'''})$$

Let $z'_{j}$ be the point in $\Omega_{j}$ (and in $\Omega_{j}'$, by (iii) ), that corresponds to the point $R_{j}\sqrt{\frac{r_{j}}{R_{j}}}+0i$ in the annulus $\{r_{j}<|z|<R_{j}\}$. By analogy we define $z''_{j}$ to be the point, that corresponds to $R_{j}\left(\frac{r_{j}}{R_{J}}\right)^{\frac{3}{10}}$.
 
Then 
$$\limsup_{\Omega\ni z\to \zeta}R_{\Omega}(z)=\limsup_{j\to\infty}R_{\Omega}(z'_{j})=2$$
and
$$\liminf_{\Omega\ni z\to \zeta}R_{\Omega}(z)=\liminf_{j\to\infty}R_{\Omega}(z''_{j})=-\infty,$$
by Theorem \ref{glownetw} and Observation \ref{biholom}.
\end{proof}
\begin{observation}$\Omega$ defined as above is hyperconvex
\end{observation}
Hyperconvexity is the same as regularity (with respect to the Dirichlet problem) in dimension $1$. One easily constructs barrier functions at each boundary point of $\Omega$ and by  Perron's method $\Omega$ is regular.

This is quite unexpected since it is known that both the Bergman kernel (see \cite{MR1210002}) and the Bergman metric (see \cite{MR1650305}, \cite{MR1714284}, \cite{MR1799743}) behave in a quite predictable way in hyperconvex domains.
  \end{section}

I would like to thank professor Zbigniew B\l ocki for numerous suggestions and encouragement.

\bibliographystyle{amsplain.bst}
\bibliography{bergmansectionalcurvature}

\end{document}